\documentclass[reqno,b5paper]{amsart}

\usepackage{amsmath}
\usepackage{amssymb}
\usepackage{amsthm}
\usepackage{enumerate}
\usepackage[mathscr]{eucal}
\usepackage{eqlist}

\theoremstyle{plain}
\newtheorem{theorem}{Theorem}[section]

\newtheorem{corollary}[theorem]{Corollary}
\theoremstyle{definition}

\newtheorem{example}[theorem]{Example}

\numberwithin{equation}{section}

\numberwithin{equation}{section}

\begin{document}
\title [Some results on singular value inequalities ] {Some results on singular value inequalities of normal operators}
\author[A. Taghavi and V. Darvish]%
{Ali Taghavi* and Vahid Darvish}

\newcommand{\acr}{\newline\indent}
\address{\llap{*\,}Department of Mathematics,\\ Faculty of Mathematical
Sciences,\\ University of Mazandaran,\\ P. O. Box 47416-1468,\\
Babolsar, Iran.} \email{taghavi@umz.ac.ir,
v.darvish@stu.umz.ac.ir }

 \subjclass[2010]{15A18, 15A42, 15A60, 47A63, 47B05}
\keywords{Singular value, compact operator, inequality, normal
operator}

\begin{abstract}
Let $x=a+ib$ be a complex number, so we have the following
inequality
$$(1/\sqrt{2})|a+b|\leq |x|\leq |a|+|b|$$
We give an operator version of above inequality. Also we obtain some
results for normal operators.
\end{abstract}
\maketitle


\section{Introduction}\label{intro}
Let $B(H)$ denote the space of all bounded linear operators on a
complex separable Hilbert space H, and let $K(H)$ denote the
two-sided ideal of compact operators in B(H). We denote the singular
values of an operator $A\in K(H)$ as $s_{1}(A)\geq s_{2}(A)\geq
\ldots$ are the eigenvalues of the positive operator
$|A|=(A^{*}A)^{1/2}$, and repeated accordingly to multiplicity. The
direct sum $A\oplus B$ denotes the block diagonal matrix
$\left[\begin{array}{cc}
A&0\\
0&B\\
\end{array}\right]$ defined on $H\oplus H$, see \cite{aud,zhan}.

For every self-adjoint operator $A$, we can write $A=A^{+}-A^{-}$
where $A^{+}$ and $A^{-}$ are positive operators. This is called the
Jordan decomposition of $A$ (see \cite{bhabook}).

 Bhatia and Kittaneh have
proved in \cite{bha} that if $A$ and $B$ are two $n\times n$
positive semidefinite matrices, then
\begin{equation}\label{bhakit}
 s_{j}(A+B)\leq \sqrt{2}s_{j}(A+iB),\ \ 1\leq j\leq n.
 \end{equation}
  Also they have shown, there exist a positive semidefinite matrix
  $A$ and Hermitian $B$ for which \eqref{bhakit} is not true.

On the other hand, Tao has proved in \cite{tao} that if $A, B, C\in
K(H)$ such that $\left[\begin{array}{cc}
A&B\\
B^{*} &C\\
\end{array}\right]\geq 0$, then
\begin{equation}\label{tao2}
2s_{j}(B)\leq s_{j}\left[\begin{array}{cc}
A&B\\
B^{*} &C\\
\end{array}\right]
\end{equation}
for $j=1,2,\ldots$

 Audeh and
Kittaneh have proved in \cite{aud}, for every $A,B,C\in K(H)$ such
that $\left[\begin{array}{cc}
A&B\\
B^{*} &C\\
\end{array}\right]\geq 0$, then
\begin{equation}\label{kit2}
s_{j}(B)\leq s_{j}(A\oplus C)
\end{equation}
for $j=1,2,\ldots$ Also by using (\ref{tao2}), they have proved that
if $A,B\in K(H)$, such that $A$ is self-adjoint, $B\geq 0$ and $\pm
A\leq B$, then
\begin{equation}\label{kit3}
2s_{j}(A)\leq s_{j}((B+A)\oplus(B-A))
\end{equation}
for $j=1,2,\ldots$

 In this paper, for each normal operator $A$,
we will prove the following inequality
\begin{equation}\label{main}
(1/\sqrt{2})s_{j}(A_{1}+A_{2})\leq s_{j}(A_{1}+iA_{2})\leq
s_{j}(|A_{1}|+|A_{2}|),
\end{equation}
for $j=1,2,\ldots$ which is an operator extension of the following
elementary and fundamental inequality
\begin{equation}\label{elec}
(1/\sqrt{2})|a+b|\leq |a+ib|\leq |a|+|b|.
\end{equation}
There exist two matrices $A$ and $B$ which the matrix inequality of
\eqref{elec} is not true.\\
As an application of \eqref{main}, we will find the upper and lower
bound for $A+iA^{*}$ for an arbitrary operator $A\in K(H)$. Also
some applications of these inequalities are obtained.

\section{Main results}
Here we will give some results for compact normal operators. For
every operator $A$, the Cartesian decomposition is to write
$A=A_{1}+iA_{2}$, where $A_{1}=\frac{A+A^{*}}{2}$ and
$A_{2}=\frac{A-A^{*}}{2i}$. If $A$ is normal operator then $A_{1}$
and $A_{2}$ commute together and vice versa.
\begin{theorem}\label{normal}
Let $A$ be a normal operator in $K(H)$. Then we have
$$(1/\sqrt{2})s_{j}(A_{1}+A_{2})\leq s_{j}(A)\leq s_{j}(|A_{1}|+|A_{2}|)$$
for $j=1,2,\ldots$ where $A=A_{1}+iA_{2}$.
\end{theorem}
\begin{proof}
We know
$\sqrt{A^{*}A}=\sqrt{(A_{1}^{2}+A_{2}^{2})+i(A_{1}A_{2}-A_{2}A_{1})}$,
since $A$ is normal we have
$\sqrt{A^{*}A}=\sqrt{A_{1}^{2}+A_{2}^{2}}$, so
$$s_{j}(A)=s_{j}(|A|)=s_{j}(\sqrt{A^{*}A})=s_{j}(\sqrt{A_{1}^{2}+A_{2}^{2}})$$
for $j=1,2,\ldots$ By using Weyl's monotonicity principle
\cite{bhabook} and the inequality
 $\sqrt{A_{1}^{2}+A_{2}^{2}}\leq
|A_{1}|+|A_{2}|$, we have the following
$$s_{j}(\sqrt{A_{1}^{2}+A_{2}^{2}})\leq s_{j}(|A_{1}|+|A_{2}|)$$
for $j=1,2,\ldots$ Now for proving left inequality, we recall this
famous inequality $0\leq (A_{1}+A_{2})^{*}(A_{1}+A_{2})\leq
2(A_{1}^{2}+A_{2}^{2})$. Therefore, by using the Weyl's monotonicity
principle we can write
$$s_{j}(\sqrt{(A_{1}+A_{2})^{*}(A_{1}+A_{2})})\leq
\sqrt{2}s_{j}(\sqrt{A_{1}^{2}+A_{2}^{2}}).$$ for $j=1,2,\ldots$  So
$$s_{j}(A_{1}+A_{2})=s_{j}(|A_{1}+A_{2}|)\leq
\sqrt{2}s_{j}(\sqrt{A_{1}^{2}+A_{2}^{2}})$$ for $j=1,2,\ldots$\qed
\end{proof}

\begin{example}\label{1}
Take $A=\left[\begin{array}{cc}
2-i&2i\\
2i&2i\\
\end{array}\right]$, by the Cartesian decomposition
we have
 $$A_{1}=\left[\begin{array}{cc}
2&0\\
0&0\\
\end{array}\right], \ \  A_{2}=\left[\begin{array}{cc}
-1&2\\
2&2\\
\end{array}\right]$$
A calculation shows, neither
$$(1/\sqrt{2})|A_{1}+A_{2}|\leq|A_{1}+iA_{1}|$$ nor
$$|A_{1}+iA_{2}|\leq |A_{1}|+|A_{2}|$$ holds.
\end{example}
\begin{example}
Let $A=\left[\begin{array}{cc}
1+i&1\\
1&i\\
\end{array}\right]$, then a calculation shows $$s_{2}(A_{1}+iA_{2})\approx
1.1756>s_{2}(|A_{1}|+|A_{2}|)\approx 0.9591.$$
\end{example}
\begin{theorem}
Let $A=A_{1}+iA_{2}$ be a normal operator in $K(H)$ such that
$-A_{2}\leq A_{1}$. Then we have
$$s_{j}(A)\leq s_{j}(2(A_{1}^{+}+A_{2}^{+})\oplus (A_{1}+A_{2}))$$
for $j=1,2,\ldots$
\end{theorem}
\begin{proof}
By Theorem \ref{normal}, we can write
$$s_{j}(A)\leq
s_{j}(2A_{1}^{+}+2A_{2}^{+}+A_{1}^{-}+A_{2}^{-}-A_{1}^{+}-A_{2}^{+})$$
for $j=1,2,\ldots$ Now imply inequality \eqref{kit3}, by supposing
$$B=2A_{1}^{+}+2A_{2}^{+}-A_{1}^{-}-A_{2}^{-}+A_{1}^{+}+A_{2}^{+}$$
hence the result follows.\qed
\end{proof}
One can easily obtain the following theorem by similar proof.
\begin{theorem}
Let $A$ be a self-adjoint operator in $K(H)$. Then
$$s_{j}(A^{+})\leq s_{j}(|A|\oplus (|A|-A)/2)$$
for $j=1,2,\ldots$  and
$$s_{j}(A^{-})\leq s_{j}(|A|\oplus (|A|+A)/2)$$
for $j=1,2,\ldots$
\end{theorem}
As mentioned in Introduction, in the following Theorem we determine
the upper and lower bound for $A+iA^{*}$.

In \cite{taghavi}, we proved the following Theorem for an arbitrary
matrix.
\begin{theorem}\label{max}
Let $A$ be a $n\times n$ complex matrix. Then
$$\sqrt{2}s_{1}(A_{1}+A_{2})\leq s_{1}(A+iA^{*})\leq 2s_{1}(A_{1}+A_{2})$$
where $A_{1}$ and $A_{2}$ are the Cartesian decomposition.
\end{theorem}
Here, we prove above Theorem for compact operator. The proof of the
following Theorem is similar to Theorem \ref{max}, but for readers
convenience we don't omit the proof.
\begin{theorem}
Let $A$ be in $K(H)$ . Then
$$\sqrt{2}s_{j}(A_{1}+A_{2})\leq s_{j}(A+iA^{*})\leq 2s_{j}(A_{1}+A_{2})$$
for $j=1,2,\ldots$ and $A=A_{1}+iA_{2}$ where $A_{1}$ and $A_{2}$
are the Cartesian decomposition.
\end{theorem}
\begin{proof}
Note that $T=A+iA^{*}$ is normal operator. On the other hand we can
write $T=T_{1}+iT_{2}$ where
$$T_{1}=((A+A^{*})+i(A^{*}-A))/2  , \quad
T_{2}=((A-A^{*})+i(A^{*}+A))/2i$$ it is enough to compare $T_{1}$
and $T_{2}$ to see $T_{1}=T_{2}$. So
\begin{equation}\label{mosavi}
T_{1}+T_{2}=(A+A^{*})+i(A^{*}-A).
\end{equation}
Now apply Theorem \ref{normal}, we have
\begin{equation}\label{111}
(1/\sqrt{2})s_{j}(T_{1}+T_{2})\leq s_{j}(T_{1}+iT_{2})\leq
s_{j}(|T_{1}|+|T_{2}|)
\end{equation}
for $j=1,2,\ldots$ Substitute (\ref{mosavi}),
$T_{1}+iT_{2}=A+iA^{*}$ and $T_{1}$ in (\ref{111}) we obtain
\begin{align*}
(1/\sqrt{2})s_{j}((A+A^{*})+i(A^{*}-A))&\leq s_{j}(A+iA^{*})\\
 &\leq
2s_{j}((A+A^{*})/2+i(A^{*}-A)/2)\\
&=s_{j}((A+A^{*})+i(A^{*}-A))\\
\end{align*}
for $j=1,2,\ldots$ By writing $A_{1}=(A+A^{*})/2$ and
$A_{2}=(A-A^{*})/2i$ we have
$$(1/\sqrt{2})s_{j}(2A_{1}+2A_{2})\leq s_{j}(A+iA^{*})\leq
s_{j}(2A_{1}+2A_{2})$$ for $j=1,2,\ldots$ Finally
$$\sqrt{2}s_{j}(A_{1}+A_{2})\leq s_{j}(A+iA^{*})\leq
2s_{j}(A_{1}+A_{2})$$ for $j=1,2,\ldots$\qed
\end{proof}

In \cite{bha1} Bhatia and Kittaneh have proved for two $n\times n$
complex matrices $A$ and $B$ we have
$$s_{j}(A^{*}B+B^{*}A)\leq
s_{j}((A^{*}A+B^{*}B)\oplus(A^{*}A+B^{*}B)), \ \ 1\leq j\leq n.$$ In
particular for Hermitian $A$ and $B$ we have the following
$$s_{j}(AB+BA)\leq
s_{j}((A^{2}+B^{2})\oplus(A^{2}+B^{2})), \ \ 1\leq j\leq n.$$ Also
Hirzallah and Kittaneh have shown in \cite{hir} that we have
$$s_{j}(AB^{*}+BA^{*})\leq
s_{j}((A^{*}A+B^{*}B)\oplus(A^{*}A+B^{*}B)), \ \ 1\leq j\leq n.$$
Here we establish some results for compact operators in $B(H)$.
\begin{theorem}
Let $A$ and $B$ be in $K(H)$. Then
$$s_{j}(AB+BA)\leq s_{j}((A^{*}A+B^{*}B)\oplus(AA^{*}+BB^{*}))$$
for $j=1,2,\ldots$
\end{theorem}
\begin{proof}
Suppose $X=\left[\begin{array}{cc}
A&B\\
B^{*}&A^{*}\\
\end{array}\right]$.
So we have,
$$XX^{*}=\left[\begin{array}{cc}
AA^{*}+BB^{*}&AB+BA\\
B^{*}A^{*}+A^{*}B^{*}&A^{*}A+B^{*}B\\
\end{array}\right]$$
By using inequality (\ref{kit2}), we obtain
$$s_{j}(AB+BA)\leq s_{j}(A^{*}A+BB^{*}\oplus AA^{*}+B^{*}B)$$
for $j=1,2,\ldots$\qed
\end{proof}
\begin{corollary}
Let $A$ and $B$ be two normal operators in $K(H)$. Then
$$s_{j}(AB+BA)\leq s_{j}(AA^{*}+BB^{*}\oplus AA^{*}+BB^{*})$$
for $j=1,2,\ldots$
\end{corollary}

In the operator theory, several extensions of the notion of the
normality are known \cite{sai}. One of the most important and most
widely studied classes among them is the hyponormality (i.e.,
$TT^{*}\leq T^{*}T$)\cite{con}. Recall that a compact hyponormal
operator in $B(H)$ is normal \cite{wat}, so we can restate these
result for compact hyponormal operators.



\begin{thebibliography}{20}
\bibitem{aud}
        {\sc W. Audeh and F. Kittaneh},
        {\it Singular value inequalities for compact operators},
       Linear Algebra Appl., {\bf 437} (2012), 2516--2522.
\bibitem{bhabook}
        {\sc R. Bhatia},
        {\it Matrix Analysis},
        Springer-Verlag, New York, 1997.



\bibitem{bha1}
        {\sc R. Bhatia and F. Kittaneh},
        {\it The matrix arithmetic-geometric mean
inequality revisited},
       Linear Algebra Appl., {\bf 428} (2008), 2177--2191.

\bibitem{bha}
        {\sc R. Bhatia and F. Kittaneh},
        {\it The singular values of $A+B$ and $A+iB$},
       Linear Algebra Appl., {\bf 431} (2009), 1502--1508.

\bibitem{con}
        {\sc J. B. Conway and W. Szymanski},
        {\it Linear combinations of hyponormal operators},
      Rocky Mountain J. Math., {\bf 18} (1988), 695--705.


\bibitem{hir}
        {\sc O. Hirzallah and F. Kittaneh},
        {\it Inequalities for sums and direct sums of Hilbert space operators},
       Linear Algebra Appl., {\bf 424} (2007), 71--82.

\bibitem{sai}
        {\sc T. Saito},
        {\it Hyponormal operators and related topics},
       Springer-Verlag, 1972.



\bibitem{taghavi}
        {\sc A. Taghavi, V. Darvish and H. M. Nazari},
        {\it Maximum singular value inequalities of
normal matrices},
        Accepted in Studies in Nonlinear Sciences.



\bibitem{tao}
        {\sc Y. Tao},
        {\it More results on singular value inequalities of matrices},
       Linear Algebra Appl., {\bf 416} (2006), 724--729.

\bibitem{wat}
        {\sc H. Watanabe},
        {\it Operators characterized by certain Cauchy-Schwarz type inequalities},
       Publ. RIMS, Kyoto Univ., {\bf 30} (1994), 249--259.


\bibitem{zhan}
        {\sc X. Zhan},
        {\it Matrix Inequalitis}, Springer-Verlag, Berlin, 2002.




\end{thebibliography}
\end{document}